# Shooting-Projection Method for Two-Point Boundary Value Problems

Stefan M. Filipov [1], Ivan D. Gospodinov [1], István Faragó [2]

[1] Department of Computer Science, Faculty of Physical, Mathematical, and Technical Sciences,
University of Chemical Technology and Metallurgy, Sofia
8 "Kl. Ohridski" Blvd., Sofia 1756, Bulgaria

[2] Department of Applied Analysis and Computational Mathematics, Faculty of Science,
Eötvös Loránd University, MTA-ELTE Research Group, Budapest
1117 Budapest Pazmany P. s. 1/C., Hungary, faragois@cs.elte.hu

**Abstract**

This paper presents a novel shooting method for solving two-point boundary value problems for second order ordinary differential equations. The method works as follows: first, a guess for the initial condition is made and an integration of the differential equation is performed to obtain an initial value problem solution; then, the end value of the solution is used in a simple iteration formula to correct the initial condition; the process is repeated until the second boundary condition is satisfied. The iteration formula is derived utilizing an auxiliary function that satisfies both boundary conditions and minimizes the $H^1$ semi-norm of the difference between itself and the initial value problem solution.

**Introduction**

Let $u(t)$ be a real-valued function of a real independent variable $t \in [a,b]$. Consider the two-point boundary value problem (TPBVP) [1], [2], [3]

$$u'' = f(t,u,u'), \ t \in (a,b),  \tag{1}$$
$$u(a) = u_a, \ u(b) = u_b,  \tag{2}$$

where $a$, $b$, $u_a$, and $u_b$ are given constants, and $f$ is a given function that specifies the differential equation (1). Let us assume that $f$ is such that the problem (1)(2) has a unique solution on the interval $[a,b]$. The basic idea of any shooting method for solving TPBVPs is to replace the boundary conditions (2) with the initial conditions

$$u(a) = u_a, \ u'(a) = V_a,  \tag{3}$$

and treat the TPBVP as an initial value problem (IVP). In (3) $V_a$ is the value of the derivative of the TPBVP solution at the first (left) boundary $t=a$. Since $V_a$ is not known, one can make a guess for its value and solve (1) together with (3), using the guess value $v_a$ instead of $V_a$. The obtained IVP solution $u(t;v_a)$ satisfies $u'(a;v_a)=v_a$ and the first boundary condition, i.e. $u(a;v_a)=u_a$, but typically it does not satisfy the second (right) boundary condition, i.e. $u(b;v_a) \neq u_b$. The difference

$$E(v_a) = u(b;v_a) - u_b  \tag{4}$$

is the *deviation* from the second boundary condition. Due to the uniqueness of the TPBVP solution, $E(v_a)=0$ if and only if $v_a=V_a$. Then, the corresponding IVP solution $u(t;V_a)$ is the sought TPBVP solution. Hence, any shooting method for TPBVPs is, in fact, a root seeking procedure for finding the root $V_a$ of $E(v_a)=0$.

Some widely used iterative methods for finding roots of algebraic equations are the secant method [4], [5], the Newton method [6], [7], the constant-slope Newton method [4], [5], and the fixed-point method [4], [8]. For all of these methods the next root approximation $v_a^*$ is found by subtracting the term $E(v_a)/k$ from the current root approximation $v_a$ (see (12) below). The value of $k$ depends on the method used. Some of the methods readjust $k$ at each iteration. For the secant or the Newton shooting methods $k$ is the slope of the current secant or tangent to $E$ line, respectively. Other methods do not readjust $k$. For the fixed-point method $k$ is a fixed, different from zero arbitrary number, usually equal to one. For the constant-slope Newton method the value of $k$ is also fixed and is equal to the slope of the tangent to $E$ line evaluated at the starting guess for $V_a$.

The proposed in this paper shooting-projection method employs an auxiliary function $u^*(t)$ which is an $H^1$-projection of the IVP solution $u(t;v_a)$. The function $u^*$ satisfies both boundary conditions and minimizes the $H^1$ semi-norm of the difference between itself and the IVP solution. As proven below, $u^*$ is an approximate TPBVP solution. It is used to derive an iteration formula for the new initial condition $v_a^*$ as a function of the old initial condition $v_a$ (eqn. 12). The formula is the same as the formula for the other considered in this work shooting methods, differing only in the value of $k$.

**Shooting-projection method**

Let the function $u(t;v_a)$, $t \in [a,b]$ be the IVP solution satisfying the initial conditions $u(a;v_a)=u_a$ and $u'(a;v_a)=v_a$, where $v_a$ is some guess for $V_a$. We transform the IVP solution $u$ into the function $u^*(t)$ ($u^* \in H^1$, $t \in [a,b]$) which satisfies the two boundary conditions (2), i.e.

$$u^*(a) = u_a, \quad u^*(b) = u_b, \tag{5}$$

and minimizes the $H^1$ semi-norm of the difference between itself and the IVP solution $u$, i.e. minimizes the functional

$$S(u^*) = \int_a^b (u^{*\prime} - u')^2 \, dt. \tag{6}$$

The function $u^*$ will be called an $H^1$-*projection* of the IVP solution $u$. If a function $u^*$ that satisfies the two boundary conditions (5) should minimize the functional (6), then the following Euler-Lagrange equation must hold:

$$\frac{d}{dt}\left(\frac{\partial (u^{*\prime}-u')^2}{\partial u^{*\prime}}\right) - \frac{\partial (u^{*\prime}-u')^2}{\partial u^*} = 0, \quad t \in (a,b). \tag{7}$$

Performing the differentiation yields

$$u^{*\prime\prime} = u'', \quad t \in (a,b). \tag{8}$$

Since the IVP solution $u$ satisfies the differential equation (1) we can replace $u''$ in (8) with $f(t, u, u')$ and then expand $f$ in Taylor series around $u^*$ to obtain

$$u^{*\prime\prime} = f(t, u^*, u^{*\prime}) + \left.\frac{\partial f}{\partial u}\right|_{u=u^*, u'=u^{*\prime}} (u - u^*) + \left.\frac{\partial f}{\partial u'}\right|_{u=u^*, u'=u^{*\prime}} (u' - u^{*\prime}) + \ldots, \quad t \in (a,b). \tag{9}$$

Equation (9) tells us that if $v_a$ is close to $V_a$, and therefore $u^*$ is close to $u$, then $u^*$ satisfies the differential equation (1) *approximately*. Since $u^*$ also satisfies the boundary conditions (5), it follows that the $H^1$-projection $u^*$ is an approximate TPBVP solution. Therefore, the derivative of $u^*$ at the first boundary could be used as a new initial condition in an iterative shooting procedure.

To obtain the value of the derivative of $u^*$ at the first boundary as a function of the initial condition $v_a$ equation (8) should be integrated. First, integrating (8) on $[a,t]$, we get

$$u^{*\prime}(t) - u^{*\prime}(a) = u'(t;v_a) - u'(a;v_a). \tag{10}$$

Then, integrating (10) on $[a,b]$, we obtain

$$u^*(b) - u^*(a) - u^{*\prime}(a)(b-a) = u(b;v_a) - u(a;v_a) - u'(a;v_a)(b-a). \tag{11}$$

Finally, introducing the notation $u^{*\prime}(a) = v_a^*$, and using the boundary conditions (5) and the initial conditions $u(a;v_a)=u_a$, $u'(a;v_a)=v_a$, expression (11) is rearranged to

$$v_a^* = v_a - \frac{E(v_a)}{k}, \tag{12}$$

where $E(v_a)$ is the deviation of the IVP solution from the second boundary condition (eqn. 4), and $k=b-a$. This is the shooting-projection iteration formula. Given an initial condition approximation $v_a$, the next approximation $v_a^*$ is obtained using formula (12). This formula is the same as the formula for the other considered in this work shooting methods, only the value of the slope $k$ is different, namely $b-a$. As in the fixed-point method and the constant-slope Newton method, the value of $k$ is fixed, however it is not arbitrary nor does it depend on the guess for the initial condition, but comes from the boundaries. As shown in the Numerical examples section, there are cases wherein using the shooting-projection method, i.e. iteration formula (12) with $k=b-a$, leads to convergence while the other two fixed-$k$ methods diverge. Compared to the methods that readjust $k$ at each iteration, such as the Newton or the secant shooting methods, the shooting-projection method may perform better in specific situations, such as the presence of inflection points or local exrema of $E(v_a)$ in the vicinity of the root [5]. Examples are provided in the Numerical examples section.

Iteration (12) is a fixed-point iteration and $v_a - E(v_a)/k$ is its associated iteration function. Sometimes, (12) is also referred to as Picard iteration or functional iteration for the solution of $E(v_a)=0$ [9]. If $|dv_a^*/dv_a|<1$ in some interval around the root of $E(v_a)=0$ and the starting initial condition $v_a$ is inside this interval, then the shooting-projection iterative procedure converges to the root. Similarly to the fixed-point method and the constant-slope Newton method, if close enough to the root the shooting-projection method converges [4], [5] when

$$|1 - m/k| < 1, \tag{13}$$

where $m$ is the slope of $E(v_a)$ at the root $V_a$. Since $k>0$, inequality (13) is satisfied if and only if $0<m<2k$. The convergence is linear except when $k=m$. Then the convergence is quadratic. A more general way to require convergence (and uniqueness of solution) is to require that the right-hand side of (12) is a contraction mapping in the whole of $\Re$ [9].

**Numerical examples**

The shooting-projection method was tested on several TPBVPs. The first example demonstrates the importance of using $k=b-a$ in formula (12) instead of the usual $k=1$ used by the fixed-point method. The second example demonstrates the stability of the proposed shooting-projection method around inflection points of $E(v_a)$. The third example demonstrates the stability of the method around local extrema of $E(v_a)$. The performance of the shooting-projection method was compared with the performance of the other considered in this work shooting methods and, for the third example, with the finite difference method (FDM) as well.

**Example 1.** Consider the TPBVP

$$u'' = \frac{1}{8}\exp(u), \ t \in (a,b), \tag{14}$$
$$u(a) = \bar{u}(a), \ u(b) = \bar{u}(b),$$

where $a = -2\sqrt{2}/3$, $b = 4\sqrt{2}/3$, and $\bar{u}(t) = \ln\left(\frac{1}{2}\tan^2\left(\frac{\sqrt{2}}{8}t + \frac{13}{6}\right) + \frac{1}{2}\right)$ is the exact solution to (14).

Starting from initial condition $v_a=0$, the TPBVP (14) was solved using the shooting-projection method. Figure 1 plots the function $E(v_a)$ (found numerically) and the sequence of parallel lines with slope $k=2\sqrt{2}$ that lead to the root of $E(v_a)=0$.

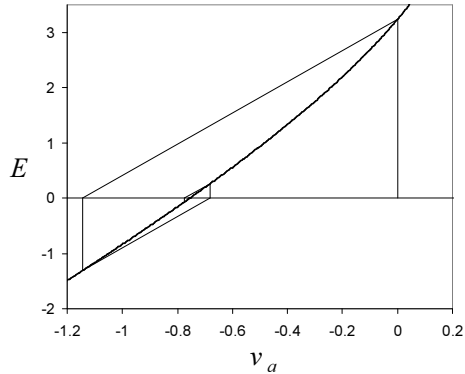

Fig. 1. Iterative solution of (14) using the shooting-projection method.

When the fixed-point method was tried with its usual slope value $k=1$ it diverged.

**Example 2.** Consider the TPBVP

$$u'' = -\frac{3u^2 u'}{t}, \ t \in (1,2), \tag{15}$$
$$u(1) = 1/\sqrt{2}, \ u(2) = 2/\sqrt{5},$$

with an exact solution $u(t) = t/\sqrt{1+t^2}$. The function $E(v_a)$ for eqn. (15) is shown in fig. 2. Due to the presence of an inflection point close to the root, the Newton and the constant-slope Newton shooting methods diverge for any initial condition $v_a$ outside a tiny interval around the root. For starting value $v_a=5$ both methods diverged in the way shown in fig. 2b and 2c, respectively. In contrast, the shooting-projection method converges for any starting initial condition $v_a$. For the same starting value $v_a=5$ it converged in 17 iterations with $|E|<0.001$ (see fig. 2a).

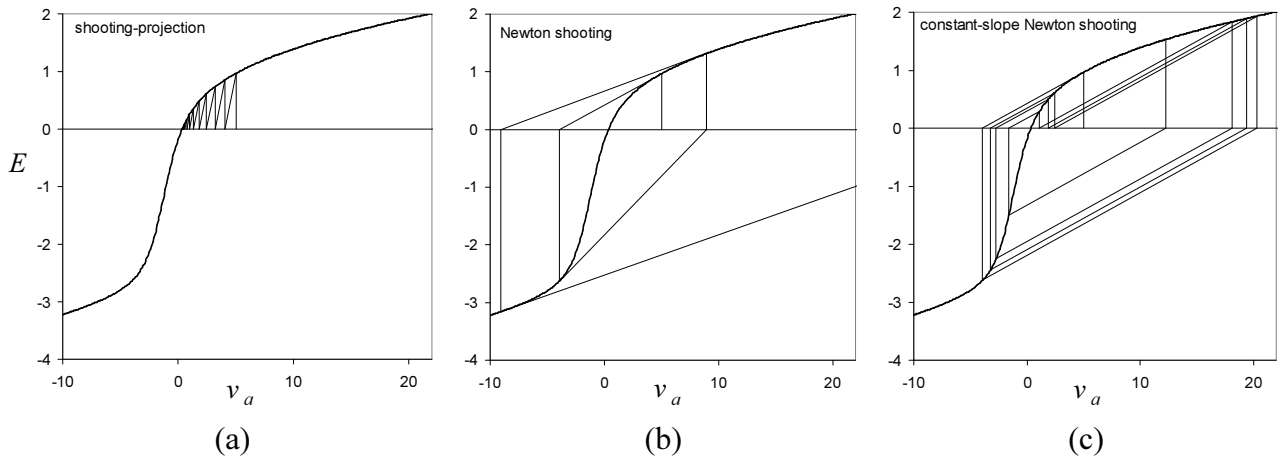

Fig. 2. Convergence/divergence sequences for the iterative solution of (15) using:
(a) the shooting-projection; (b) the Newton; and (c) the constant-slope Newton shooting method.

**Example 3**. Consider the TPBVP

$$u'' = -\frac{1}{50} u \cosh\left(\frac{tu}{5} + u\right), \; t \in (0,5), \tag{16}$$
$$u(0) = 1, \; u(5) = 2.$$

Starting from $v_a=0$ the shooting-projection method finds a solution in 14 iterations with $|E|<0.0001$. The IVP solution at each iteration is shown in fig. 3.

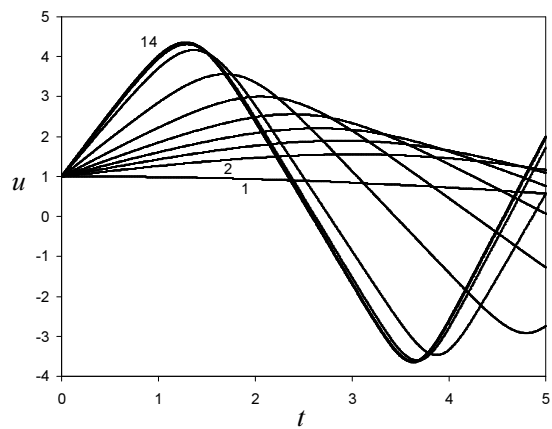

Fig. 3. Evolution of the IVP solution when solving (16) using the shooting-projection method. The TPBVP solution is reached in 14 iterations.

Figure 4a shows the sequence of parallel lines with slope $k=5$ that lead to the root of $E(v_a)=0$. The obtained root value is $V_a=3.2232$. The consecutive values of $v_a$ at each iteration, found using formula (12), are indicated with 1,2,3,… in the figure. They correspond to the IVP solutions shown in fig.3. Figure 4a shows that the function $E(v_a)$ has two local extrema between the starting guess $v_a=0$ and the root. Unlike other methods, the shooting-projection method does not have convergence issues around such extrema. It is easy to see that the method will find the root for any starting $v_a$ inside the interval shown in the figure.

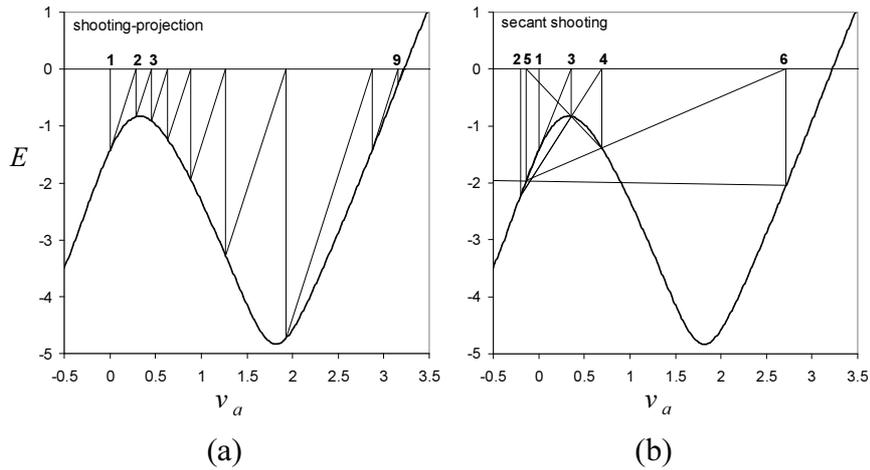

Fig. 4. Convergence/divergence sequences for the iterative solution of (16) using:
(a) the shooting-projection method; and (b) the secant method.

The secant method started from the secant line through points $(-0.2, E(-0.2))$ and $(0, E(0))$ diverges in the way shown in fig. 4b. The Newton shooting method started from $v_a=0$ also fails to find the root. The fixed-point method with $k=1$ diverges. The constant-slope Newton method, if started from any $v_a$ corresponding to a negative slope of $E$ or to a positive slope of $E$ less than $m/2$ [4], [5] will also fail to find the root.

To solve the TPBVP (16) the FDM was also tried. The differential equation was discretized and the obtained algebraic equations were complemented with the two boundary conditions. To solve the resulting system the Jacobi, Gauss-Seidel, Picard, and Newton iterative schemes were applied. All four schemes failed.

**Conclusion**

This paper described a novel shooting method for iterative solution of TPBVPs. An $H^1$-projection of the IVP solution was used to derive an iteration formula (eqn. 12) for correcting the initial condition. The main feature of the iteration formula is that the slope $k$ is fixed with a value that comes from the boundaries and thus it is neither arbitrary nor does it depend on the choice of initial conditions. It was demonstrated that there are cases wherein the proposed shooting-projection method finds a solution while the other considered in this work shooting and finite difference methods do not.